\documentstyle [epic,eepic,11pt]{amsart}

\newtheorem{theorem}{Theorem}[section]                                    
\newtheorem{Theorem}[theorem]{Theorem}                                   
\newtheorem{definition}[theorem]{Definition}
 
\newtheorem{lemma}[theorem]{Lemma}

\newtheorem{Example}[theorem]{Example}  

\newcommand{\pr} {\smallskip\noindent{\em Proof.\,\,}}

\newenvironment{proof}	{\pr}{\hspace*{\fill}\qed\\}

\newcommand\R{{\Bbb R}}

\newcommand\C{{\Bbb C}}
\newcommand\CP{{\Bbb C}{\Bbb P}}

\newcommand\Z{{\Bbb Z}}

\newcommand\X{{\cal X}}
\newcommand\Y{{\cal Y}}

\newcommand\fh {{\frak h}}
\newcommand\fk {{\frak k}}
\newcommand\ft {{\frak t}}

\begin{document}

\title{Examples of non-K\"ahler \\ Hamiltonian torus actions}
\author{Susan Tolman}\date{\today}
\address{Department of Mathematics, Mass.\ Inst.\ of Tech., Cambridge, MA 
02139} 
\email{tolman@@math.mit.edu}
\thanks{The author was partially supported by an NSF
postdoctoral fellowship}
\maketitle

\begin{abstract}
An important question with a rich history is the extent to which the
symplectic category is larger than the K\"ahler category.
Many interesting examples of non-K\"ahler symplectic manifolds have
been constructed \cite{T} \cite{M} \cite{G}.
However, sufficiently large symmetries can force
a symplectic manifolds to be K\"ahler \cite{D} \cite{K}.
In this paper, we solve several outstanding problems by
constructing the first symplectic manifold with large non-trivial symmetries
which does not admit an invariant K\"ahler structure.
The proof that it is not K\"ahler is based on the Atiyah-Guillemin-Sternberg
convexity theorem \cite{At} \cite{GS}.
Using the ideas of this paper, C. Woodward shows that even the symplectic
analogue of spherical varieties need not be K\"ahler \cite{W}.
\end{abstract}

\section{Introduction}

A {\bf symplectic manifold} is a manifold $M$ with a nondegenerate closed
two form $\omega$.  
A complex structure $J : TM \to TM$ is {\bf compatible} with $\omega$ 
if  $g(X,Y) := \omega(J(X),Y)$ is a Reimanian metric.
A {\bf K\"ahler structure} 
is a symplectic form $\omega$ and a compatible complex structure $J$ on $M$.
Many theorems are much easier to prove for K\"ahler manifolds
than for symplectic manifolds. Indeed, some are only known to be true in 
the K\"ahler case, and others are known to be true only in the K\"ahler case.
Therefore, an important question with a rich history 
is the extent to which the symplectic category is larger 
than the K\"ahler category.

Many interesting examples of non-K\"ahler symplectic manifolds have been 
constructed.
For example, in 1971 Thurston
constructed a compact four dimensional symplectic manifold 
which is not K\"ahler \cite{T}. 
In 1984, D.\ McDuff used Thurston's example and symplectic blowups
to construct a compact ten dimensional {\em simply connected}
non-K\"ahler symplectic manifold \cite{M}.
There are two ways to show that these manifolds do not admit K\"ahler
structures;
they each have an odd Betti number which is odd,
and  their cohomology does not have the hard Lefshetz property.
In 1994, 
among many other examples,   
R. Gompf constructed a symplectic form on a Gompf-Mrowka manifold,
a four dimensional simply connected symplectic manifold
which does not admit any complex structure \cite{G} \cite{GM}.

Let $T$ be a torus which acts on a manifold $M$, preserving a symplectic
form $\omega$.  The action is  {\bf Hamiltonian} if there
exists a {\bf moment map} 
$\phi: M \to \ft^*$ such that $i_{\xi_M}\omega = - d \left< \phi, \xi
\right>$ for all $\xi \in \ft$, where $\xi_M$ is the induced vector field on
$M$.  
For mathematicians who study Hamiltonian actions, the important
question is to what extent do these symmetries force symplectic manifolds
to admit invariant K\"ahler structures.

On the one hand,  it is easy to find a
non-K\"ahler symplectic manifolds with a Hamiltonian group action. 
let $M$  be a symplectic manifold which does not admit a K\"ahler structure.
The action of $S^1$ on $M \times \CP^1$ given by
$ \lambda \cdot (m, [x, y]) = (m, [\lambda x, y])$ is
Hamiltonion, and $M'$ does not admit an invariant K\"ahler structure.
However, this example is not very interesting.
The components of the fixed point set and the reduced spaces
are themselves symplectic manifolds which do not admit K\"ahler structures.  

More interestingly, 
E.\ Lerman used McDuff's example to construct a compact twelve
dimensional simply connected symplectic manifold and a Hamiltonian
circle action on it with an isolated fixed point \cite{Lt}. 
The manifold is not K\"ahler because 
it has an odd dimensional odd Betti number. 
In this example, some of the fixed point sets and some of the reduced
spaces are themselves non-K\"ahler.

On the other hand,
sufficiently large symmetries can force
a symplectic manfiods to be K\"ahler.
let a torus $T$ act effectively on 
a compact symplectic manifold $(M,\omega)$ in a Hamiltonian fashion.
T.\ Delzant shows that if the dimension of $T$ is half 
the dimension of $M$, then $(M,\omega)$  admits a $T$ invariant  compatible 
complex structure \cite{D}.
Similarly, if $M$ is four dimensional and $T$ is one dimensional,
then Y.\ Karshon shows that
$(M,\omega)$ admits a $T$ invariant compatible complex structure \cite{Ka}.
(See also \cite{Au} and \cite{AH}).

In this paper, 
we give the first symplectic manifold with large non-trivial symmetries
which does not admit an invariant K\"ahler structure.
We form a six dimensional compact symplectic  
manifold with an effective Hamiltonian two dimensional torus action. 
Using Atiyah-Guillemin-Sternberg convexity theorem \cite{At} \cite{GS},
we develop a new criterion for showing that a 
symplectic manifold does not admit an compatible invariant complex structure,  
and use this to show that our example does not admit any invariant 
K\"ahler structure.
Moreover, using these techniques it is possible to find many other
symplectic manifolds with similar properties.

This example solves several outstanding existence problems.
For instance, it shows that a symplectic manifold with a Hamiltonian torus
action need not admit an invariant K\"ahler structure even
if all the reduced spaces are K\"ahler.
It is the first known example with isolated fixed points.
This example is sharp in the sense that it is not possible to find either
an example with a lower dimensional symplectic manifold,
or an example with a larger Hamiltonian torus action (relative to the
symplectic manifold).
Additionally, this manifold is simply connected, and the action is quasi-free.

However, the most important implication of this paper is that
the category of symplectic manifolds with
Hamiltonian torus actions is much richer than the category 
of K\"ahler manifolds with compatible Hamiltonian torus actions. 
In recent years their has been a degree of convergence between
the study of symplectic manifolds with Hamiltonian group actions,
and the study of algebraic varieties with algebraic group actions.
Some mathematicians had feared -- or hoped --
that traditional Hamiltonian techniques might be rendered obsolete.
After all, more powerful techniques are often available in algebraic
geometry.  However, I believe that this family of examples shows that
K\"ahler manifolds, far from being the rule, are merely an
interesting special case -- even among highly symmetric symplectic manifolds.

Using the ideas of this paper, C. Woodward shows that even the symplectic
analogue of spherical varieties need not be K\"ahler \cite{W}.
He  finds an alternate construction for my example
and for several others, using symplectic cuts.
In so doing, he shows that they in fact admit
a multiplicity free Hamiltonian $U(2)$ action.
F.\ Knop has independently constructed other symplectic manifolds
which admit a multiplicity free Hamiltonian $G$ action, where
$G$ is a  compact nonabelian group.
He shows that they do not admit compatible $G$ 
invariant complex structures by using the 
classification of K\"ahler
manifolds which admit multiplicity free  Hamiltonian actions \cite{K}. 

In Section \ref{s_xray}, we  define the 
x-ray of 
a compact symplectic manifold with
a Hamiltonian torus action. 

In Section \ref{s_criterion}, 
we prove that if the x-ray of a symplectic manifold
with a Hamiltonian torus action does not satisfy the extension criterion,
then the symplectic manifold does not admit a compatible invariant complex
structure.  

In Section \ref{s_example}, 
we glue together pieces of two K\"{a}hler manifolds 
to construct a symplectic manifold with a Hamiltonian torus action.
Because its x-ray does not satisfy the extension criterion,
it does not admit an invariant compatible complex structure.

In Section \ref{s_varying}, we show that the manifold which we constructed
in Section \ref{s_example} does not admit any $T$ invariant K\"ahler structure.  

In Section \ref{s_more},  we describe
how to construct more symplectic manifolds with Hamiltonian torus actions
which do not admit compatible invariant complex structures.
\subsection{Acknowledgements}

I am happy to have this opportunity to thank Yael Karshon for all that
she taught me about symplectic geometry.  In particular, my ideas about
x-rays, two torus actions on six manifolds, and symplectic toric manifolds
all grew out of joint projects with her, and made this paper possible.

Similarly, I would like to thank Chris Woodward for his enthusiasm and
for his invaluable help with many aspects of this project .

I am grateful to Eugene Lerman for introducing me to this problem and
for offering me mathematical and practical advice. 

\section{ The x-ray}
\label{s_xray}

In this section, we  
define the  x-ray of  a compact symplectic manifold with a Hamiltonian
torus action.
In the rest of the paper, we will use the x-ray both to construct an example,
and to show that it does not admit an invariant K\"ahler structure.

Let a torus $T$ act on a compact symplectic manifold $M$ with  moment
map $\phi: M \to \ft^*$.
For each subgroup $K \subset T$, 
let  $M^K$ be the set of points fixed by $K$,  
and let $\X_K$ be the set of connected components of $M^K$.
The {\bf closed orbit type stratification} of $M$ is the set $\X = \cup_{K \subset T} \X_K$;
this set is  partial ordered by inclusion.

By equivariant Darboux, 
every $X \in \X$ is itself is a symplectic $T$ invariant manifold
with moment map $\phi|_X$.
Therefore, by the Atiyah-Guillemin-Sternberg 
convexity theorem \cite{At} \cite{GS},
$\phi(X) \subset \ft^*$ is a convex polytope.

\begin{definition}
Let a torus $T$ act on a compact symplectic manifold $(M,\omega)$ with  moment
map $\phi : M \to \ft^*$.
The {\bf x-ray} of $(M,\omega,\phi)$ is the closed orbit type stratification
$\X$ of $M$, and the convex polytope $\phi(X)$ for each $X \in \X$.
\end{definition}

Notice that we do not need to understand the geometry of each $X$;
$\X$ simply functions as a (partially ordered) index set.
More abstractly, {\bf an x-ray} is  any 
partially ordered set  with a convex polytope associated to each element.

\begin{Example}
\label{cp2can}{\em
Define a manifold
$$ \tilde{M} = \left\{ [x_0,x_1,x_2,y_0,y_1,y_2] \in \CP^5 \left| x_i y_j^4 = x_j y_i^4
\mbox{ for all } i \mbox{ and } j \right. \right\}.$$ 
Let $\tilde{\omega}$ be the Fubini-Studi symplectic form.
The torus $S^1 \times S^1$ acts on $\tilde{M}$ by 
$$(\alpha,\beta) \cdot [x_0,x_1,x_2,y_0,y_1,y_2] =
 [ x_0,\alpha^4 x_1, \beta^4 x_2, \alpha \beta y_0, \alpha^2 \beta y_1, \alpha
 \beta^2y_2].$$
The x-ray for a  moment map $\tilde{\phi}: M \to \R^2$ appears in the left
hand side of Figure~1.\footnote
{Indeed,  a three dimensional torus acts effectively on $\tilde{M}$, that is
$\tilde{M}$ is a symplectic toric manifold.
The action described above is the action of a two dimensional subtorus.
The above x-ray is thus the projection of a three dimensional polytope.}

The small grey dots denote the weight lattice; the
dot in the bottom left is the origin.
Each large black dot is the image of a 
connected component of the fixed point set.
Each black line 
is the image of a connected component of the set
of points with a given one dimensional stabilizer.
The dashed line will be explained later.
}\end{Example}

\begin{Example}
\label{cp12} {\em
Let $S^1 \times S^1$ act on $\hat{M} = \CP^1 \times \CP^2$
by $$(\alpha,\beta) \cdot \left([x_0, x_1], [y_0, y_1, y_2]\right) = 
\left([\alpha x_0, x_1], [\alpha y_0,  \beta y_1, y_2]\right).$$
Let $\hat{\omega} = \pi_1^*(\omega_1) + 3*\pi_2^*(\omega_2)$,
where $\pi_i$ is the projection map onto the $i$'th componet, and
$\omega_i$ is the Fubini-Studi symplectic form on $\CP^i$.
The x-ray for a  moment map $\hat{\phi}: M \to \R^2$ appears in the right
hand side of Figure~1.\footnote
{The footnote above also applies to $\hat{M}$.}
}\end{Example}

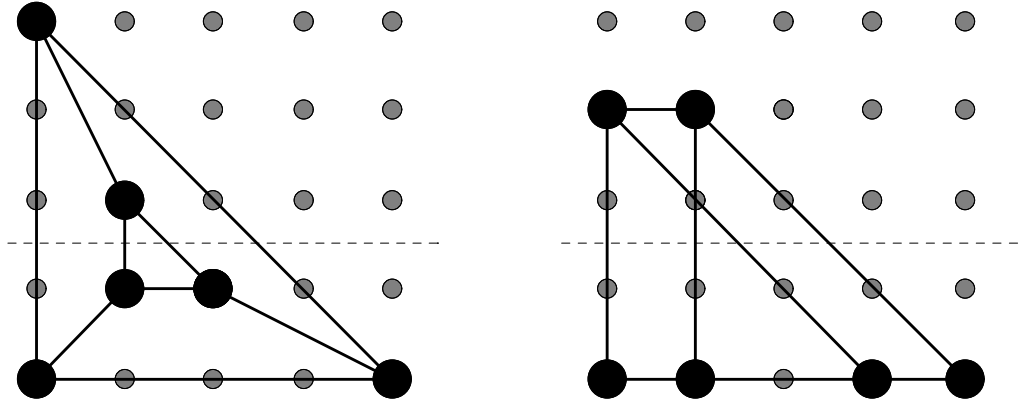
\begin{figure}
\setlength{\unitlength}{0.01250000in}
\begingroup\makeatletter\ifx\SetFigFont\undefined%
\gdef\SetFigFont#1#2#3#4#5{%
  \reset@font\fontsize{#1}{#2pt}%
  \fontfamily{#3}\fontseries{#4}\fontshape{#5}%
  \selectfont}%
\fi\endgroup%
{\renewcommand{\dashlinestretch}{30}
\begin{picture}(427,179)(0,-10)
\texture{44555555 55aaaaaa aa555555 55aaaaaa aa555555 55aaaaaa aa555555 55aaaaaa 
	aa555555 55aaaaaa aa555555 55aaaaaa aa555555 55aaaaaa aa555555 55aaaaaa 
	aa555555 55aaaaaa aa555555 55aaaaaa aa555555 55aaaaaa aa555555 55aaaaaa 
	aa555555 55aaaaaa aa555555 55aaaaaa aa555555 55aaaaaa aa555555 55aaaaaa }
\put(14,45){\shade\ellipse{8}{8}}
\put(14,45){\ellipse{8}{8}}
\put(14,82){\shade\ellipse{8}{8}}
\put(14,82){\ellipse{8}{8}}
\put(14,120){\shade\ellipse{8}{8}}
\put(14,120){\ellipse{8}{8}}
\put(51,157){\shade\ellipse{8}{8}}
\put(51,157){\ellipse{8}{8}}
\put(51,120){\shade\ellipse{8}{8}}
\put(51,120){\ellipse{8}{8}}
\put(51,82){\shade\ellipse{8}{8}}
\put(51,82){\ellipse{8}{8}}
\put(51,45){\shade\ellipse{8}{8}}
\put(51,45){\ellipse{8}{8}}
\put(88,7){\shade\ellipse{8}{8}}
\put(88,7){\ellipse{8}{8}}
\put(88,45){\shade\ellipse{8}{8}}
\put(88,45){\ellipse{8}{8}}
\put(88,82){\shade\ellipse{8}{8}}
\put(88,82){\ellipse{8}{8}}
\put(88,120){\shade\ellipse{8}{8}}
\put(88,120){\ellipse{8}{8}}
\put(88,157){\shade\ellipse{8}{8}}
\put(88,157){\ellipse{8}{8}}
\put(126,157){\shade\ellipse{8}{8}}
\put(126,157){\ellipse{8}{8}}
\put(126,120){\shade\ellipse{8}{8}}
\put(126,120){\ellipse{8}{8}}
\put(126,82){\shade\ellipse{8}{8}}
\put(126,82){\ellipse{8}{8}}
\put(126,45){\shade\ellipse{8}{8}}
\put(126,45){\ellipse{8}{8}}
\put(126,7){\shade\ellipse{8}{8}}
\put(126,7){\ellipse{8}{8}}
\put(163,7){\shade\ellipse{8}{8}}
\put(163,7){\ellipse{8}{8}}
\put(163,45){\shade\ellipse{8}{8}}
\put(163,45){\ellipse{8}{8}}
\put(163,82){\shade\ellipse{8}{8}}
\put(163,82){\ellipse{8}{8}}
\put(163,120){\shade\ellipse{8}{8}}
\put(163,120){\ellipse{8}{8}}
\put(163,157){\shade\ellipse{8}{8}}
\put(163,157){\ellipse{8}{8}}
\put(163,7){\blacken\ellipse{16}{16}}
\put(163,7){\ellipse{16}{16}}
\put(88,45){\blacken\ellipse{16}{16}}
\put(88,45){\ellipse{16}{16}}
\put(88,45){\blacken\ellipse{16}{16}}
\put(88,45){\ellipse{16}{16}}
\put(51,45){\blacken\ellipse{16}{16}}
\put(51,45){\ellipse{16}{16}}
\put(51,82){\blacken\ellipse{16}{16}}
\put(51,82){\ellipse{16}{16}}
\put(253,7){\shade\ellipse{8}{8}}
\put(253,7){\ellipse{8}{8}}
\put(253,45){\shade\ellipse{8}{8}}
\put(253,45){\ellipse{8}{8}}
\put(253,82){\shade\ellipse{8}{8}}
\put(253,82){\ellipse{8}{8}}
\put(253,120){\shade\ellipse{8}{8}}
\put(253,120){\ellipse{8}{8}}
\put(253,157){\shade\ellipse{8}{8}}
\put(253,157){\ellipse{8}{8}}
\put(290,157){\shade\ellipse{8}{8}}
\put(290,157){\ellipse{8}{8}}
\put(290,120){\shade\ellipse{8}{8}}
\put(290,120){\ellipse{8}{8}}
\put(290,82){\shade\ellipse{8}{8}}
\put(290,82){\ellipse{8}{8}}
\put(290,45){\shade\ellipse{8}{8}}
\put(290,45){\ellipse{8}{8}}
\put(290,7){\shade\ellipse{8}{8}}
\put(290,7){\ellipse{8}{8}}
\put(327,7){\shade\ellipse{8}{8}}
\put(327,7){\ellipse{8}{8}}
\put(327,45){\shade\ellipse{8}{8}}
\put(327,45){\ellipse{8}{8}}
\put(290,82){\shade\ellipse{8}{8}}
\put(290,82){\ellipse{8}{8}}
\put(327,82){\shade\ellipse{8}{8}}
\put(327,82){\ellipse{8}{8}}
\put(327,120){\shade\ellipse{8}{8}}
\put(327,120){\ellipse{8}{8}}
\put(327,120){\shade\ellipse{8}{8}}
\put(327,120){\ellipse{8}{8}}
\put(327,157){\shade\ellipse{8}{8}}
\put(327,157){\ellipse{8}{8}}
\put(364,157){\shade\ellipse{8}{8}}
\put(364,157){\ellipse{8}{8}}
\put(364,120){\shade\ellipse{8}{8}}
\put(364,120){\ellipse{8}{8}}
\put(364,82){\shade\ellipse{8}{8}}
\put(364,82){\ellipse{8}{8}}
\put(364,45){\shade\ellipse{8}{8}}
\put(364,45){\ellipse{8}{8}}
\put(364,7){\shade\ellipse{8}{8}}
\put(364,7){\ellipse{8}{8}}
\put(403,7){\shade\ellipse{8}{8}}
\put(403,7){\ellipse{8}{8}}
\put(403,45){\shade\ellipse{8}{8}}
\put(403,45){\ellipse{8}{8}}
\put(403,82){\shade\ellipse{8}{8}}
\put(403,82){\ellipse{8}{8}}
\put(403,120){\shade\ellipse{8}{8}}
\put(403,120){\ellipse{8}{8}}
\put(403,157){\shade\ellipse{8}{8}}
\put(403,157){\ellipse{8}{8}}
\put(253,120){\blacken\ellipse{16}{16}}
\put(253,120){\ellipse{16}{16}}
\put(290,120){\blacken\ellipse{16}{16}}
\put(290,120){\ellipse{16}{16}}
\put(253,7){\blacken\ellipse{16}{16}}
\put(253,7){\ellipse{16}{16}}
\put(290,7){\blacken\ellipse{16}{16}}
\put(290,7){\ellipse{16}{16}}
\put(403,7){\blacken\ellipse{16}{16}}
\put(403,7){\ellipse{16}{16}}
\put(364,7){\blacken\ellipse{16}{16}}
\put(364,7){\ellipse{16}{16}}
\put(14,157){\blacken\ellipse{16}{16}}
\put(14,157){\ellipse{16}{16}}
\put(51,7){\shade\ellipse{8}{8}}
\put(51,7){\ellipse{8}{8}}
\put(14,7){\blacken\ellipse{16}{16}}
\put(14,7){\ellipse{16}{16}}
\dashline{4.000}(234,64)(425,64)
\Thicklines
\path(14,157)(163,7)
\path(14,157)(14,7)(163,7)
	(88,45)(51,45)(51,82)(14,157)
\path(88,45)(51,82)
\path(14,7)(51,45)
\path(253,120)(253,7)(403,7)
	(290,120)(253,120)(364,7)
\path(290,7)(290,120)
\thinlines
\dashline{4.000}(2,64)(182,64)(182,64)
\end{picture}
}
\caption{ The x-rays for Example~\ref{cp2can} and Example~\ref{cp12}}
\end{figure}

\section{The extension criterion}
\label{s_criterion}

In this section, we define the etension criterion.
We then  prove that if the x-ray of a symplectic manifold
with a Hamiltonian torus action does not satisfy this criterion,
then the symplectic manifold does not admit a compatible invariant complex
structure.

A convex polytope $\Delta \subset \ft^*$ is {\bf compatible} with an x-ray
$(\X,\phi)$ if 
\begin{enumerate}
\item
for each face $\sigma$ of $\Delta $,  there exists
$X_\sigma \in \X$ such that $\dim(\phi(X_\sigma)) = \dim(\sigma)$
and $\sigma \subseteq \phi(X_\sigma)$;
\item
if $\sigma$ and $\sigma'$ are faces of $\Delta$ such that
$\sigma \subset \sigma' $, then 
$X_\sigma \subset  X_{\sigma'}$.  
\end{enumerate}

Let a torus $T$ act on a symplectic manifold $(M,\omega)$ with moment map
$\phi: M \to \ft^*$.
If  $N \subset M$ is a symplectic $T$-invariant submanifold,
then the convex polytope $\phi(N)$ is compatible with the x-ray of $(M,\omega,\phi)$. 
There are exactly five two dimensional (and nine one dimensional)
convex polytopes which are 
compatible with the x-rays in  Examples \ref{cp2can} and \ref{cp12}.

Similarly, a convex cone $C \subset \ft^*$ is {\bf compatible} with an x-ray $\X$ if
there is a neighborhood $U$ of the vertex of $C$ such that 
\begin{enumerate}
\item
for each face $\sigma$ of $C$,  there exists
$X_\sigma \in \X$ such that $\dim(\phi(X_\sigma)) = \dim(\sigma)$
and $\sigma \cap U \subseteq \phi(X_\sigma)$, 
and 

\item
if $\sigma$ and $\sigma'$ are faces of $C$ such that
$\sigma \subset \sigma'$, then 
$X_\sigma \subset X_{\sigma'}$.  
\end{enumerate}

A convex polytope $\Delta$ is an {\bf extension} of a
convex cone $C$ if there is a neighborhood $U$ of its vertex 
such that $C \cap U = \Delta \cap U$.

\begin{definition}
\label{definition of extension criterion}
An x-ray satisfies the {\bf extension criterion} if
every compatible strictly 
convex cone can be extended to a compatible convex polytope.
\end{definition}

The x-rays in Figure~1 satisfy the extension criterion.

\begin{Example} {\em
\label{nonkahler} 
The x-ray in Figure~2 does not satisfy the extension criterion.
For instance, the cone 
$$\{ (s,t) \in \R^2 \mid  t \geq 1  \mbox{ and } s +t \leq 3 \}$$
is compatible with the x-ray, but it does not extend to a compatible
polytope.  (Remember that the origin is in the bottom left corner).
Notice, however, that below the dashed line this x-ray resembles the 
x-ray on the left in Figure~1, and that above the dashed line this
x-ray resembles the x-ray on the right in Figure~1.
}\end{Example}

\begin{figure}
\setlength{\unitlength}{0.00083333in}
\begingroup\makeatletter\ifx\SetFigFont\undefined%
\gdef\SetFigFont#1#2#3#4#5{%
  \reset@font\fontsize{#1}{#2pt}%
  \fontfamily{#3}\fontseries{#4}\fontshape{#5}%
  \selectfont}%
\fi\endgroup%
\renewcommand{\dashlinestretch}{30}
%\begin{picture}(2803,2473)(0,-10)
\begin{picture}(2803,2473)(-1800,-10)
\texture{44000000 aaaaaa aa000000 8a888a 88000000 aaaaaa aa000000 888888 
	88000000 aaaaaa aa000000 8a8a8a 8a000000 aaaaaa aa000000 888888 
	88000000 aaaaaa aa000000 8a888a 88000000 aaaaaa aa000000 888888 
	88000000 aaaaaa aa000000 8a8a8a 8a000000 aaaaaa aa000000 888888 }
%\shade\dashline{100}(1,689)(1381,689)(16,2054)
%	(1,689)
%\drawline(1,689)(1381,689)(16,2054)
%	(1,689)
\texture{44555555 55aaaaaa aa555555 55aaaaaa aa555555 55aaaaaa aa555555 55aaaaaa 
	aa555555 55aaaaaa aa555555 55aaaaaa aa555555 55aaaaaa aa555555 55aaaaaa 
	aa555555 55aaaaaa aa555555 55aaaaaa aa555555 55aaaaaa aa555555 55aaaaaa 
	aa555555 55aaaaaa aa555555 55aaaaaa aa555555 55aaaaaa aa555555 55aaaaaa }
\put(226,119){\shade\ellipse{120}{120}}
\put(226,119){\ellipse{120}{120}}
\put(226,689){\shade\ellipse{120}{120}}
\put(226,689){\ellipse{120}{120}}
\put(226,1259){\shade\ellipse{120}{120}}
\put(226,1259){\ellipse{120}{120}}
\put(226,1829){\shade\ellipse{120}{120}}
\put(226,1829){\ellipse{120}{120}}
\put(796,2399){\shade\ellipse{120}{120}}
\put(796,2399){\ellipse{120}{120}}
\put(796,1829){\shade\ellipse{120}{120}}
\put(796,1829){\ellipse{120}{120}}
\put(796,1259){\shade\ellipse{120}{120}}
\put(796,1259){\ellipse{120}{120}}
\put(796,689){\shade\ellipse{120}{120}}
\put(796,689){\ellipse{120}{120}}
\put(796,119){\shade\ellipse{120}{120}}
\put(796,119){\ellipse{120}{120}}
\put(1366,119){\shade\ellipse{120}{120}}
\put(1366,119){\ellipse{120}{120}}
\put(1366,689){\shade\ellipse{120}{120}}
\put(1366,689){\ellipse{120}{120}}
\put(1366,1259){\shade\ellipse{120}{120}}
\put(1366,1259){\ellipse{120}{120}}
\put(1366,1829){\shade\ellipse{120}{120}}
\put(1366,1829){\ellipse{120}{120}}
\put(1366,2399){\shade\ellipse{120}{120}}
\put(1366,2399){\ellipse{120}{120}}
\put(1921,2399){\shade\ellipse{120}{120}}
\put(1921,2399){\ellipse{120}{120}}
\put(1921,1829){\shade\ellipse{120}{120}}
\put(1921,1829){\ellipse{120}{120}}
\put(1921,1259){\shade\ellipse{120}{120}}
\put(1921,1259){\ellipse{120}{120}}
\put(1921,689){\shade\ellipse{120}{120}}
\put(1921,689){\ellipse{120}{120}}
\put(1921,119){\shade\ellipse{120}{120}}
\put(1921,119){\ellipse{120}{120}}
\put(2491,119){\shade\ellipse{120}{120}}
\put(2491,119){\ellipse{120}{120}}
\put(2491,689){\shade\ellipse{120}{120}}
\put(2491,689){\ellipse{120}{120}}
\put(2491,1259){\shade\ellipse{120}{120}}
\put(2491,1259){\ellipse{120}{120}}
\put(2491,1829){\shade\ellipse{120}{120}}
\put(2491,1829){\ellipse{120}{120}}
\put(2491,2399){\shade\ellipse{120}{120}}
\put(2491,2399){\ellipse{120}{120}}
\put(2491,119){\blacken\ellipse{240}{240}}
\put(2491,119){\ellipse{240}{240}}
\put(1366,689){\blacken\ellipse{240}{240}}
\put(1366,689){\ellipse{240}{240}}
\put(1366,689){\blacken\ellipse{240}{240}}
\put(1366,689){\ellipse{240}{240}}
\put(796,689){\blacken\ellipse{240}{240}}
\put(796,689){\ellipse{240}{240}}
\put(226,119){\blacken\ellipse{240}{240}}
\put(226,119){\ellipse{240}{240}}
\put(226,2399){\shade\ellipse{120}{120}}
\put(226,2399){\ellipse{120}{120}}
\put(226,1829){\blacken\ellipse{240}{240}}
\put(226,1829){\ellipse{240}{240}}
\put(796,1829){\blacken\ellipse{240}{240}}
\put(796,1829){\ellipse{240}{240}}
%\drawline(16,659)(16,659)
\dashline{60.000}(16,959)(2791,959)
%\drawline(226,1829)(226,1829)
\Thicklines
\path(796,1829)(226,1829)(226,119)
	(2491,119)(796,1829)(796,689)(226,119)
\path(2491,119)(1366,689)
\path(1366,689)(796,689)
\path(1366,689)(226,1829)
\end{picture}
\caption{ An x-ray which does not satisfy the extension criterion }  
\end{figure}
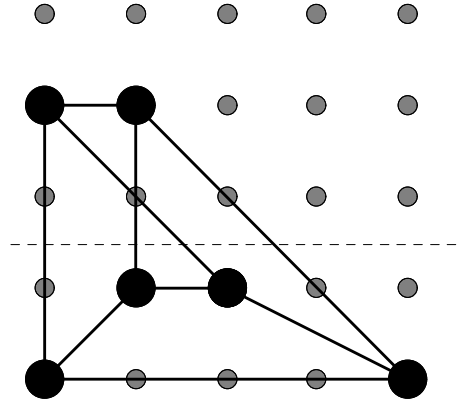

\begin{Theorem}
\label{theorem}
Let a torus $T$ act on a compact symplectic manifold $(M,\omega)$
with moment map $\phi: M \to \ft^*$.
If the x-ray of $(M,\omega,\phi)$ does not satisfy the extension criterion
(see Definition~\ref{definition of extension criterion}),
then $(M,\omega)$ does not admit a compatible $T$-invariant complex structure.
\end{Theorem}

\begin{proof}
This theorem follows directly from Theorem~\ref{Atiyah}, 
which is a reformulation of Theorem~2 in \cite{At},
and Lemma \ref{lemorbit}. 
\end{proof}

The converse need not be true.

\begin{Theorem}
\label{Atiyah}
{\em {\bf (Atiyah)}} Let a torus $T$ act by holomorphic
symplectomorphisms on 
a compact K\"ahler manifold $(M,\omega,J)$
with  moment map $\phi: M \to \ft^*$.
Then $T_\C$, the complexification of $T$, also acts on $M$.
For $y \in M$, 
the set $ \phi(\overline{T_\C \cdot y})$ is a convex polytope which 
is compatible with the x-ray of $(M,\omega,\phi)$,
where $\overline{T_\C \cdot y}$ is the closure of the $T_\C$ orbit of $y$. 
\end{Theorem}

\noindent {\em Remark.}\
In the case we need, the action of $T$ on $M$ is quasi-free.
Therefore,  $\overline{T_\C \cdot y}$ is a smooth symplectic manifold,
and $\phi(\overline{T_\C \cdot y})$ is a convex polytope by
the Atiyah-Guillemin-Sternberg convexity theorem \cite{At} \cite{GS}.

\begin{lemma}
\label{lemorbit}
Let a torus $T$ act by holomorphic symplectomorphisms  on a compact 
K\"ahler manifold $(M,\omega,J)$ with  moment map $\phi: M \to \ft^*$.
Let $C$ be a strictly convex cone which is compatible with the x-ray
of $(M,\omega,\phi)$.
There exist  $y \in M$ such that
$\phi(\overline{T_\C \cdot y})$ is an extension of $C$,
where $\overline{T_\C \cdot y}$ is the closure of the $T_\C$ orbit of $y$, 
and $T_\C$ is the complexification of $T$.
\end{lemma}

\begin{proof}
Since $C$ is compatible with the x-ray of $(M,\omega,\phi)$, there
is a fixed point component $F \subset M$ which corresponds to 
the vertex of $C$.
Given $m \in F$,
let  $\xi_1,\ldots,\xi_n \in \ft^*$ be the weights for the action of
$T_\C$ on $T_mM$. 
Let $T_\C$ act on $\C^n$ by
$t \cdot (z_1,\ldots,z_n) = (t^{\xi_1} z_1 ,\ldots,t^{\xi_n} z_n)$.
There exists a $T_\C$ equivariant holomorphic diffeormorphism $\psi$
from a neighborhood of $0$ in $\C^n$ to a neighborhood of $m$ in $M$.
(see \cite{S}).

Define $\omega':= \psi^*(\omega)$  and  $\phi' := \psi^*(\phi)$. 
Given $J \subset \{1,\ldots,n\}$, 
let $E_J = \{ z \in \C^n \mid z_i \neq 0 \mbox{ if and only if }  i \in J \}$. 
The image $C_J := \phi'(\overline{E_J})$ is equal to 
$\{ \sum_{j \in J} a_j \xi_j \mid a_j \geq 0 \}$.
Since $C$ is compatible with the x-ray of $(M,\omega,\phi)$
it is also compatible with the x-ray of $(\C^n,\omega',\phi')$.
Therefore, $C = C_J$ for some unique minimal $J$.  

Choose $y' \in  E_J$.
On the one hand,
$\phi'(\overline{T_\C \cdot y'})  \subset \phi'(\overline{E_J}) = C_J$. 
On the other hand, for every $j \in  J$, $C_{\{j\}}$ is a proper face of
$C_J$. Therefore, there exists $\eta \in \ft$
such that $\left< \eta, \xi_j \right> = 0$ but $\left< \eta, \xi_i \right> < 0$
for every  $i \in J$ such that $i \neq j$. 
It follows that  $E_{\{j\}} \subset \overline{T_\C \cdot y'}$,
and hence $C_{\{j\}} \subset \phi'(\overline{T_\C \cdot y'})$,
for all $j \in J$.
Since $\phi'(\overline{T_\C \cdot y'})$  is convex,
$\phi'(\overline{T_\C \cdot y'}) =   C_J = C$. 

Finally,  define  $y := \psi^{-1}(y')$. 
Since ${\phi'}^{-1}(\phi'(0)) = 0$,
there exists a neighborhood $U$ of the vertex of $C$
such that $\phi_{T_\C \cdot y} \cap U = C \cap U$.
\end{proof}

\section{Constructing an example}  
\label{s_example}

In this section, 
we glue together pieces of two K\"{a}hler manifolds 
to construct a symplectic manifold with a Hamiltonian torus action.
Because its x-ray does not satisfy the extension criterion,
it does not admit an invariant compatible complex structure.

\begin{lemma}
\label{lemma exists}
There exists a 
compact six dimensional manifold $(M,\omega)$
and a two dimensional torus $T$ which acts effectively
in a Hamiltonian fashion on $M$  such that $(M,\omega)$ does not admit 
any compatible invariant complex structure.   
Additionally, 
\begin{enumerate}
\item
the manifold $M$ is simply connected,
\item
the stabilizer of every point is connected,
\item
all the fixed points are isolated, and
\item
all the reduced spaces are symplectomorphic to $\CP^1$.
In particular, they all admit compatible complex structures.
\end{enumerate}
\end{lemma}

\begin{proof}
The basic idea of the construction is very simple.
Define $(\tilde{M},\tilde{\omega},\tilde{\phi})$
and $(\hat{M},\hat{\omega},\hat{\phi})$ 
as in Example \ref{cp2can} and Example \ref{cp12}, respectively.
We show that there exists a neighborhood $W \subset \R^2$ 
of the dashed line in these
examples such that 
$\tilde\phi^{-1}(W)$ and 
$\hat\phi^{-1}(W)$ are isomorphic. 
We define a new symplectic manifold by
gluing the ``bottom'' of $\tilde{M}$
to the ``top'' of $\hat{M}$.\footnote{Alternatively, this manifold can
be constructed using Gompf gluing.}

The first step is to show that if we 
reduce both spaces along the dashed line  $D = \{(s,t) \in \R^2 \mid t = 1.5 \}$,
the resulting spaces are isomorphic.
More precisely, define  $H := e \times S^1 \subset S^1 \times S^1$.
The moment map for the $H$ action is the projection of
the $S^1 \times S^1$ moment map onto its second component.
Moreover, $H$ acts freely on $\tilde{\phi}^{-1}(D)$ and $\hat\phi^{-1}(D)$.
Therefore, the reduced spaces  
$\tilde{M}//H = \tilde{\phi}^{-1}(D)/H$ and 
$\hat{M}//H = \hat{\phi}^{-1}(D)/H$ 
are smooth symplectic manifolds on which $T/H$ acts effectively in a Hamiltonian
fashion. 

The $T/H$ action on both $\tilde{M}//H$ and $\hat{M}//H$,
is free except for at four isolated fixed points. 
Moreover, the image of these four points under
the $T/H$ moment map is the same for both reduced spaces.
Applying the  classification of four dimensional
dimensional manifolds with a Hamiltonian circle action 
given by Y. Karshon,
the reduced spaces
are equivariantly symplectomorphic \cite{Ka}. (See also \cite{Au} and
\cite{AH}). 

Alternatively, $\tilde{M}$ and $\hat{M}$ are symplectic toric manifolds with
associated three dimensional polytopes.
The reduced spaces are the symplectic toric manifolds which are
associated with the cross section of these polytopes.
Specifically,
$\tilde{M}//H$ is equivariantly symplectomorphic to $\CP^1 \times \CP^1$ with the
diagonal action  
and  $\hat{M}//H$ is equivariantly symplectomorphic to 
the Hirzebruch surface $\{ [x,y,w,z] \in \CP^3  \mid x y^3 = w z^3 \}$ 
with the $S^1$ action given by $\lambda \cdot  [x,y,w,z] 
= [ x, \lambda^3 y, \lambda w, \lambda^2 z]$.
(The symplectic form is not the Fubini-Studi form in either case.)
It is straightforward to write down explicitly an $S^1$
equivariant symplectomorphism between these two spaces\footnote
{In contrast, these two spaces are not holomorphically equivalent.
Therefore, the complex structures on $\tilde{M}$ and $\hat{M}$ will not
glue together to create  a complex structure on $M$.}.

Next, we show that there exists an equivariant diffeomorphism from 
$\tilde{\phi}^{-1}(D)$ to $\hat{\phi}^{-1}(D)$ 
such that $\tilde\omega$ is the pull back of $\hat\omega$.
Because $\tilde{M}//H$ and $\hat{M}//H$ are equivariantly
symplectomorphic, it is enough to show that
$\tilde\phi^{-1}(D)$ and 
$\hat\phi^{-1}(D)$  are isomorphic as $T/H$ equivariant line bundles.

For any equivariant line bundle 
over $\tilde{M}//H = \hat{M}//H = \CP^1 \times \CP^1$,
there is a representation of $T/H$ on the fiber of each fixed point in the base
manifold. The Chern class of the line bundle 
is determined once we know the representation.
(For example, the integral of the Chern class over the first 
sphere is given by the difference of the representations at $(N,N)$ and $(S,N)$,
where $S$ and $N$ denote the north and south poles, respectively.)

Moreover, these fixed-point representations are determined 
by their kernels.  Their kernels, in turn, are determined by the stabilizers
of points in $\tilde{\phi}^{-1}(D)$ and $\hat{\phi}^{-1}(D)$.
Since all stabilizers are connected, they can be read from the x-ray.
Since $\tilde{M}$ and $\hat{M}$ have identical x-rays  near $D$,
they have the same representations at corresponding fixed points,
therefore there exists an equivariant diffeomorphism from 
$\tilde{\phi}^{-1}(D)$ to $\hat{\phi}^{-1}(D)$ 
such that $\tilde\omega$ is the pull back of $\hat\omega$.

By the equivariant coisotropic embedding theorem 
there exists neighborhoods of
$\tilde\phi^{-1}(D)$ and $\tilde{\phi}^{-1}(D)$ which are
equivariantly symplectomorphic.
Because $\tilde{M}$ and $\hat{M}$ are compact
there exists a neighborhood $W$ of $D$ such that 
there exists a
equivariant symplectomorphism
$f: \tilde\phi^{-1}(W) \to \hat{\phi}^{-1}(W)$.

Define open sets  $U := W \cup  \{(s,t) \in \R^2 \mid t < 1.5 \}$ and
$V := W \cup  \{(s,t) \in \R^2 \mid t > 1.5 \}$.
We define a new compact symplectic manifold  $(M,\omega)$ with 
an effective Hamiltonian $T$ action by 
$$ M := \tilde\phi^{-1}(U) \coprod \hat{\phi}^{-1}(V) / \sim \ ,$$
where $x \sim f(x)$ for all $x \in \tilde\phi^{-1}(W)$.

Over $U$, the x-rays of $M$ and $\tilde{M}$ agree;
over $V$, the x-rays of $M$ and $\hat{M}$ agree.
The x-ray of $(M,\omega,\phi)$ is 
the one drawn in Example \ref{nonkahler}, which does
not satisfy the extension criterion.
By Theorem~\ref{theorem},
$(M,\omega)$  does not admit a compatible $T$ invariant complex
structure.  

To see that $M$ is simply connected, it is enough to check 
that the $H$ moment map is
a Morse function and that the index of every critical point is even.

The other claims follow from the construction.
\end{proof}

\section{Varying the symplectic form}
\label{s_varying}

In this section, we show that the manifold which we constructed
in Section \ref{s_example}
does not admit any $T$-invariant K\"ahler structure.  

Let a torus $T$ act effectively on a compact manifold $M$.
There may be many invariant symplectic forms $\omega$ on $M$, 
each with its own moment map $\phi:M \to \ft^*$.
However, the x-ray of $(M,\omega,\phi)$ 
only depends in a limited way on the symplectic form $\omega$.
The closed orbit type stratification
$\X$ depends {\em only} on the action of $T$ on $M$.
In contrast, 
the polytope $\phi(X) \subset \ft^*$ for $X \in \X$ 
does depend on the cohomology class of
$\omega$ (though not on $\omega$ itself).  
By the Atiyah-Guillemin-Sternberg convexity theorem \cite{At} \cite{GS}, 
$\phi(X)$ is the convex hull of $\{ \phi(F) \mid F \subset X
\hbox{ and } F \in \X_{\{1\}} \}$, 
where $\X_{\{1\}}$ denotes the set of components of
the fixed point set.
Therefore, the x-ray is determined by 
$\phi(F)$ for $F \in \X_{\{1\}}$. 
In fact, given $F$ and $F'$ in $X_{\{1\}}$ and 
$X \in \X_K$ such that $F \subset X$ and $F' \subset X$,
$\phi(F)$ and $\phi(F')$ must satisfy
$\phi(F) - \phi(F') \in \fk^\circ \subset \ft^* $, where $\fk^\circ$ denotes
the annihilator of the lie algebra of $K$.
Therefore, the x-ray is determined by even less information.

\begin{theorem}
There exists a 
compact six dimensional symplectic manifold $(M,\omega)$
and a two dimensional torus $T$ which acts effectively on $M$
in a Hamiltonian fashion such that $M$ does not admit 
any invariant K\"ahler structure.   
Additionally, 
\begin{enumerate}
\item
the manifold $M$ is simply connected,
\item
the stabilizer of every point is connected,
\item
all the fixed points are isolated, and
\item
all the reduced spaces are symplectomorphic to $\CP^1$.
In particular, they all admit K\"ahler structures.
\end{enumerate}
\end{theorem}

\begin{proof}
Let $(M,\omega,\phi)$ be the symplectic manifold constructed in 
Lemma~\ref{lemma exists}.
Let $F$ and $F'$ denote the components of the fixed point set
such that $\phi(F) = (0,0)$ and $\phi(F') = (1,3)$.

Now choose {\em any} $S^1 \times S^1$ invariant symplectic form $\omega'$ on $M$.
Since $M$ is simply connected, there is a moment map $\phi': M \to \R^2$.
To prove the theorem, we must show that the x-ray of $(M,\omega',\phi')$
fails to satisfy the extension criterion.

Applying the discussion in the previous paragraph,
this x-ray is determined by $\phi'(F)$ and $\phi'(F')$.
Moreover, since the extension criterion is not affected by translations, 
we may assume that $\phi'(F)$ is zero.  
Therefore, the x-ray is determined by $\phi'(F') = (s,t)$.
In particular, the components of the fixed point set which are mapped by
$\phi$ to $(0,0), (4,0), (1,1), (2,1), (0,3)$, and $(1,3)$, 
will be mapped by $\phi'$ 
to $(0,0), (s+t,0), (s,s) (t-s,s), (0,t)$ and $(s,t)$, respectively.
We may further restrict attention to the case $s \geq 0$ because
the x-rays for $\phi'(F') = (s,t)$ and $\phi'(F') = (-s,-t)$ are
mirror images.   

Notice that the extension criterion does not depend on the length  of the edges;
it can only change as one of the edges ``switches direction''.
One or more of the edges will have zero length exactly if
$s = 0, t = 0, t= -s, t=s $, or $t = 2s$.
This leaves us with five cases:
$0 < t < s,\ s < t < 2s,\ 0 < t < s,\ -s < t < 0,$ and $t < -s$.
The first case has already been discussed.  Representative samples of
the other cases appear in Figure~3; $\phi'(F) = (0,0)$
is denoted by a grey dot in a large black dot.
Note that $0 < t < s$ is impossible because
the overall shape is not convex.  The other x-rays may be possible but 
do not satisfy the extension property.
\end{proof}

\begin{figure}
\setlength{\unitlength}{0.01in}
%\begin{picture}(542,466)(30,-10)
\begin{picture}(542,421)(30,30)
\texture{cccccccc 0 0 0 cccccccc 0 0 0 
	cccccccc 0 0 0 cccccccc 0 0 0 
	cccccccc 0 0 0 cccccccc 0 0 0 
	cccccccc 0 0 0 cccccccc 0 0 0 }
\put(171,160){\shade\ellipse{10}{10}}
\put(171,160){\ellipse{10}{10}}
\put(171,120){\shade\ellipse{10}{10}}
\put(171,120){\ellipse{10}{10}}
\put(131,120){\shade\ellipse{10}{10}}
\put(131,120){\ellipse{10}{10}}
\put(131,160){\shade\ellipse{10}{10}}
\put(131,160){\ellipse{10}{10}}
\put(91,160){\shade\ellipse{10}{10}}
\put(91,160){\ellipse{10}{10}}
\put(91,120){\shade\ellipse{10}{10}}
\put(91,120){\ellipse{10}{10}}
\put(51,120){\shade\ellipse{10}{10}}
\put(51,120){\ellipse{10}{10}}
\put(51,160){\shade\ellipse{10}{10}}
\put(51,160){\ellipse{10}{10}}
\put(11,200){\shade\ellipse{10}{10}}
\put(11,200){\ellipse{10}{10}}
\put(11,160){\shade\ellipse{10}{10}}
\put(11,160){\ellipse{10}{10}}
\put(11,120){\shade\ellipse{10}{10}}
\put(11,120){\ellipse{10}{10}}
\put(11,80){\shade\ellipse{10}{10}}
\put(11,80){\ellipse{10}{10}}
\put(51,80){\shade\ellipse{10}{10}}
\put(51,80){\ellipse{10}{10}}
\put(91,80){\shade\ellipse{10}{10}}
\put(91,80){\ellipse{10}{10}}
\put(131,80){\shade\ellipse{10}{10}}
\put(131,80){\ellipse{10}{10}}
\put(171,80){\shade\ellipse{10}{10}}
\put(171,80){\ellipse{10}{10}}
\put(211,80){\shade\ellipse{10}{10}}
\put(211,80){\ellipse{10}{10}}
\put(211,120){\shade\ellipse{10}{10}}
\put(211,120){\ellipse{10}{10}}
\put(211,160){\shade\ellipse{10}{10}}
\put(211,160){\ellipse{10}{10}}
\put(211,200){\shade\ellipse{10}{10}}
\put(211,200){\ellipse{10}{10}}
\put(171,200){\shade\ellipse{10}{10}}
\put(171,200){\ellipse{10}{10}}
\put(131,200){\shade\ellipse{10}{10}}
\put(131,200){\ellipse{10}{10}}
\put(91,200){\shade\ellipse{10}{10}}
\put(91,200){\ellipse{10}{10}}
\put(51,200){\shade\ellipse{10}{10}}
\put(51,200){\ellipse{10}{10}}
\put(211,200){\blacken\ellipse{22}{22}}
\put(211,200){\ellipse{22}{22}}
\put(211,80){\blacken\ellipse{22}{22}}
\put(211,80){\ellipse{22}{22}}
\put(171,120){\blacken\ellipse{22}{22}}
\put(171,120){\ellipse{22}{22}}
\put(131,120){\blacken\ellipse{22}{22}}
\put(131,120){\ellipse{22}{22}}
\put(131,80){\blacken\ellipse{22}{22}}
\put(131,80){\ellipse{22}{22}}
\put(11,200){\blacken\ellipse{22}{22}}
\put(11,200){\ellipse{22}{22}}
\put(491,160){\shade\ellipse{10}{10}}
\put(491,160){\ellipse{10}{10}}
\put(491,120){\shade\ellipse{10}{10}}
\put(491,120){\ellipse{10}{10}}
\put(451,120){\shade\ellipse{10}{10}}
\put(451,120){\ellipse{10}{10}}
\put(451,160){\shade\ellipse{10}{10}}
\put(451,160){\ellipse{10}{10}}
\put(411,160){\shade\ellipse{10}{10}}
\put(411,160){\ellipse{10}{10}}
\put(411,120){\shade\ellipse{10}{10}}
\put(411,120){\ellipse{10}{10}}
\put(371,120){\shade\ellipse{10}{10}}
\put(371,120){\ellipse{10}{10}}
\put(371,160){\shade\ellipse{10}{10}}
\put(371,160){\ellipse{10}{10}}
\put(331,200){\shade\ellipse{10}{10}}
\put(331,200){\ellipse{10}{10}}
\put(331,160){\shade\ellipse{10}{10}}
\put(331,160){\ellipse{10}{10}}
\put(331,120){\shade\ellipse{10}{10}}
\put(331,120){\ellipse{10}{10}}
\put(331,80){\shade\ellipse{10}{10}}
\put(331,80){\ellipse{10}{10}}
\put(371,80){\shade\ellipse{10}{10}}
\put(371,80){\ellipse{10}{10}}
\put(411,80){\shade\ellipse{10}{10}}
\put(411,80){\ellipse{10}{10}}
\put(451,80){\shade\ellipse{10}{10}}
\put(451,80){\ellipse{10}{10}}
\put(491,80){\shade\ellipse{10}{10}}
\put(491,80){\ellipse{10}{10}}
\put(531,80){\shade\ellipse{10}{10}}
\put(531,80){\ellipse{10}{10}}
\put(531,120){\shade\ellipse{10}{10}}
\put(531,120){\ellipse{10}{10}}
\put(531,160){\shade\ellipse{10}{10}}
\put(531,160){\ellipse{10}{10}}
\put(531,200){\shade\ellipse{10}{10}}
\put(531,200){\ellipse{10}{10}}
\put(491,200){\shade\ellipse{10}{10}}
\put(491,200){\ellipse{10}{10}}
\put(451,200){\shade\ellipse{10}{10}}
\put(451,200){\ellipse{10}{10}}
\put(411,200){\shade\ellipse{10}{10}}
\put(411,200){\ellipse{10}{10}}
\put(371,200){\shade\ellipse{10}{10}}
\put(371,200){\ellipse{10}{10}}
\put(531,200){\blacken\ellipse{22}{22}}
\put(531,200){\ellipse{22}{22}}
\put(491,160){\blacken\ellipse{22}{22}}
\put(491,160){\ellipse{22}{22}}
\put(531,80){\blacken\ellipse{22}{22}}
\put(531,80){\ellipse{22}{22}}
\put(491,80){\blacken\ellipse{22}{22}}
\put(491,80){\ellipse{22}{22}}
\put(451,160){\blacken\ellipse{22}{22}}
\put(451,160){\ellipse{22}{22}}
\put(371,200){\blacken\ellipse{22}{22}}
\put(371,200){\ellipse{22}{22}}
\put(491,400){\shade\ellipse{10}{10}}
\put(491,400){\ellipse{10}{10}}
\put(491,360){\shade\ellipse{10}{10}}
\put(491,360){\ellipse{10}{10}}
\put(451,360){\shade\ellipse{10}{10}}
\put(451,360){\ellipse{10}{10}}
\put(451,400){\shade\ellipse{10}{10}}
\put(451,400){\ellipse{10}{10}}
\put(411,400){\shade\ellipse{10}{10}}
\put(411,400){\ellipse{10}{10}}
\put(411,360){\shade\ellipse{10}{10}}
\put(411,360){\ellipse{10}{10}}
\put(371,360){\shade\ellipse{10}{10}}
\put(371,360){\ellipse{10}{10}}
\put(371,400){\shade\ellipse{10}{10}}
\put(371,400){\ellipse{10}{10}}
\put(331,440){\shade\ellipse{10}{10}}
\put(331,440){\ellipse{10}{10}}
\put(331,400){\shade\ellipse{10}{10}}
\put(331,400){\ellipse{10}{10}}
\put(331,360){\shade\ellipse{10}{10}}
\put(331,360){\ellipse{10}{10}}
\put(331,320){\shade\ellipse{10}{10}}
\put(331,320){\ellipse{10}{10}}
\put(371,320){\shade\ellipse{10}{10}}
\put(371,320){\ellipse{10}{10}}
\put(411,320){\shade\ellipse{10}{10}}
\put(411,320){\ellipse{10}{10}}
\put(451,320){\shade\ellipse{10}{10}}
\put(451,320){\ellipse{10}{10}}
\put(491,320){\shade\ellipse{10}{10}}
\put(491,320){\ellipse{10}{10}}
\put(531,320){\shade\ellipse{10}{10}}
\put(531,320){\ellipse{10}{10}}
\put(531,360){\shade\ellipse{10}{10}}
\put(531,360){\ellipse{10}{10}}
\put(531,400){\shade\ellipse{10}{10}}
\put(531,400){\ellipse{10}{10}}
\put(531,440){\shade\ellipse{10}{10}}
\put(531,440){\ellipse{10}{10}}
\put(491,440){\shade\ellipse{10}{10}}
\put(491,440){\ellipse{10}{10}}
\put(451,440){\shade\ellipse{10}{10}}
\put(451,440){\ellipse{10}{10}}
\put(411,440){\shade\ellipse{10}{10}}
\put(411,440){\ellipse{10}{10}}
\put(371,440){\shade\ellipse{10}{10}}
\put(371,440){\ellipse{10}{10}}
\put(531,320){\blacken\ellipse{22}{22}}
\put(531,320){\ellipse{22}{22}}
\put(371,320){\blacken\ellipse{22}{22}}
\put(371,320){\ellipse{22}{22}}
\put(371,360){\blacken\ellipse{22}{22}}
\put(371,360){\ellipse{22}{22}}
\put(331,400){\blacken\ellipse{22}{22}}
\put(331,400){\ellipse{22}{22}}
\put(491,360){\blacken\ellipse{22}{22}}
\put(491,360){\ellipse{22}{22}}
\put(491,400){\blacken\ellipse{22}{22}}
\put(491,400){\ellipse{22}{22}}
\put(171,400){\shade\ellipse{10}{10}}
\put(171,400){\ellipse{10}{10}}
\put(171,360){\shade\ellipse{10}{10}}
\put(171,360){\ellipse{10}{10}}
\put(131,360){\shade\ellipse{10}{10}}
\put(131,360){\ellipse{10}{10}}
\put(131,400){\shade\ellipse{10}{10}}
\put(131,400){\ellipse{10}{10}}
\put(91,400){\shade\ellipse{10}{10}}
\put(91,400){\ellipse{10}{10}}
\put(91,360){\shade\ellipse{10}{10}}
\put(91,360){\ellipse{10}{10}}
\put(51,360){\shade\ellipse{10}{10}}
\put(51,360){\ellipse{10}{10}}
\put(51,400){\shade\ellipse{10}{10}}
\put(51,400){\ellipse{10}{10}}
\put(11,440){\shade\ellipse{10}{10}}
\put(11,440){\ellipse{10}{10}}
\put(11,400){\shade\ellipse{10}{10}}
\put(11,400){\ellipse{10}{10}}
\put(11,360){\shade\ellipse{10}{10}}
\put(11,360){\ellipse{10}{10}}
\put(11,320){\shade\ellipse{10}{10}}
\put(11,320){\ellipse{10}{10}}
\put(51,320){\shade\ellipse{10}{10}}
\put(51,320){\ellipse{10}{10}}
\put(91,320){\shade\ellipse{10}{10}}
\put(91,320){\ellipse{10}{10}}
\put(131,320){\shade\ellipse{10}{10}}
\put(131,320){\ellipse{10}{10}}
\put(171,320){\shade\ellipse{10}{10}}
\put(171,320){\ellipse{10}{10}}
\put(211,320){\shade\ellipse{10}{10}}
\put(211,320){\ellipse{10}{10}}
\put(211,360){\shade\ellipse{10}{10}}
\put(211,360){\ellipse{10}{10}}
\put(211,400){\shade\ellipse{10}{10}}
\put(211,400){\ellipse{10}{10}}
\put(211,440){\shade\ellipse{10}{10}}
\put(211,440){\ellipse{10}{10}}
\put(171,440){\shade\ellipse{10}{10}}
\put(171,440){\ellipse{10}{10}}
\put(131,440){\shade\ellipse{10}{10}}
\put(131,440){\ellipse{10}{10}}
\put(91,440){\shade\ellipse{10}{10}}
\put(91,440){\ellipse{10}{10}}
\put(51,440){\shade\ellipse{10}{10}}
\put(51,440){\ellipse{10}{10}}
\put(11,440){\blacken\ellipse{22}{22}}
\put(11,440){\ellipse{22}{22}}
\put(91,440){\blacken\ellipse{22}{22}}
\put(91,440){\ellipse{22}{22}}
\put(91,400){\blacken\ellipse{22}{22}}
\put(91,400){\ellipse{22}{22}}
\put(51,400){\blacken\ellipse{22}{22}}
\put(51,400){\ellipse{22}{22}}
\put(11,320){\blacken\ellipse{22}{22}}
\put(11,320){\ellipse{22}{22}}
\put(211,320){\blacken\ellipse{22}{22}}
\put(211,320){\ellipse{22}{22}}
\put(11,320){\shade\ellipse{10}{10}}
\put(11,320){\ellipse{10}{10}}
\put(371,320){\shade\ellipse{10}{10}}
\put(371,320){\ellipse{10}{10}}
\put(131,120){\shade\ellipse{10}{10}}
\put(131,120){\ellipse{10}{10}}
\put(491,160){\shade\ellipse{10}{10}}
\put(491,160){\ellipse{10}{10}}
\Thicklines
\path(211,320)(11,320)(11,440)
	(91,440)(211,320)(51,400)
	(91,400)(11,320)
\path(11,440)(51,400)
\path(91,440)(91,400)
\path(11,200)(211,200)(211,80)
	(131,80)(11,200)(171,120)
	(131,120)(211,200)
\path(131,120)(131,80)
\path(171,120)(211,80)
\path(371,320)(371,360)(331,400)
	(491,400)(491,360)(531,320)
	(371,320)(491,400)
\path(331,400)(531,320)
\path(371,360)(491,360)
\path(531,200)(531,80)(491,80)
	(371,200)(531,200)(491,160)
	(451,160)(371,200)
\drawline(491,160)(491,160)
\path(491,160)(491,80)
\path(451,160)(531,80)
\put(91,280){\makebox(0,0)[lb]{ $s < t < 2s$ }}
\put(91,40){\makebox(0,0)[lb]{ $ -s < t < 0$ }}
\put(411,280){\makebox(0,0)[lb]{ $0 < t < s$ }}
\put(411,40){\makebox(0,0)[lb]{ $ t < -s$ }}
\end{picture}
\caption{Possible x-rays for other symplectic forms on $M$}
\end{figure}
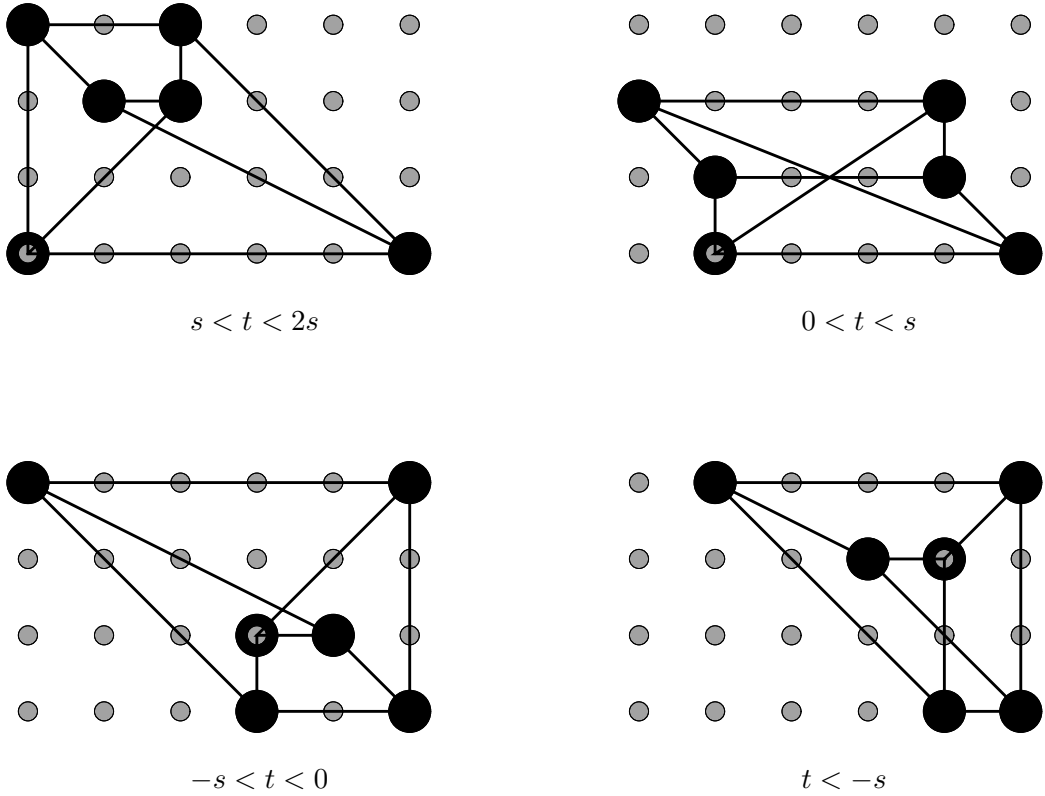

\noindent {\em Remark.}\ 
Because every odd Betti number is zero,
we cannot prove that $M$ is not K\"ahler using that route.
Additionally, as C.\ Woodward has pointed out,  $M$ does satisfy the 
hard Lefshetz property.
At the time of this writing, it is not know whether or not $M$ admits
a  K\"ahler structure.

\section{More examples}
\label{s_more}
In this section,  we describe
how to construct many more symplectic manifolds with Hamiltonian torus actions
which do not admit compatible invariant complex structures.
These profuse examples show that the category of symplectic manifolds with
Hamiltonian torus actions is much richer than the category 
of K\"ahler manifolds with compatible Hamiltonian torus actions.

The construction in Section~\ref{s_example} can be extended
to create other symplectic manifolds with Hamiltonian
torus actions which do not admit compatible invariant complex
structures.
Instead of beginning with Example \ref{cp2can} and \ref{cp12},
let a two dimensional torus $T$
act in a Hamiltonian fashion on any two six dimensional symplectic manifolds,
which we'll denote by $(\tilde{M},\tilde{\omega})$ and $(\hat{M},\hat{\omega})$.
Consider a subtorus  $H \subset T$ and $\eta \in \ft^*$. 
Define $D := \{ \alpha \in \ft^* | \pi(\alpha) = \eta \}$,
where $\pi : \ft^* \to \fh^*$ is the projection map.
The construction in Lemma~\ref{lemma exists} can be applied to
create a new manifold  when the following conditions are satisfied: 
The reduction of  $\tilde{M}$ by $H$ at $\eta$ and the
reduction of $\hat{M}$ by $H$ at $\eta$ are smooth symplectic
manifolds. 
The circle $T/H$ acts on both reduced spaces with isolated fixed points.
The x-rays of 
$(\tilde{M},\tilde{\omega},\tilde{\phi})$ and 
$(\hat{M},\hat{\omega},\hat{\phi})$ are identical near $D$,
and the corresponding elements of the closed orbit type stratifications
have the same stabilizer subgroups.
With a little patience, 
it is possible to construct an abundance of  manifolds this way.   
Typically, they do not satisfy the extension criterion.
Moreover, the arguments in Section~\ref{s_varying} can be adapted
to show that many of them do not admit any invariant K\"ahler structure.

Another way to create examples is to start 
with a symplectic manifold $(M,\omega)$ with a Hamiltonian torus action
which does not admit an invariant K\"ahler structure.
The product of $(M,\omega)$ and any
symplectic manifold on which a torus acts with isolated fixed points
will have $M$ itself as a closed orbit type strata.
Therefore, it cannot admit an invariant K\"ahler structure.

\end{document}